\documentclass{kluwer}
\begin{document}
\begin{article}
\begin{opening}
\title{Asymptotic analysis of the GI/M/1/n loss system as
n increases to infinity}

\author{\surname{Vyacheslav M. Abramov}\email{vyachesl@inter.net.il}}
\institute{24/6 Balfour st., Petach Tiqva 49350, Israel}




\runningtitle{Asymptotic analysis of the GI/M/1/n loss system}
\runningauthor{VYACHESLAV ABRAMOV}



\begin{abstract}
This paper provides the asymptotic analysis of the loss probability
in the $GI/M/1/n$ queueing
system as $n$ increases to infinity.
The approach of this paper is alternative to that of the recent
papers of Choi and Kim [8] and Choi et al [9] and based on application of modern
Tauberian theorems with remainder. This enables us to simplify
the proofs of the results on asymptotic behavior of the loss probability
of the abovementioned paper of Choi and Kim [9] as well as to obtain some new results.
\end{abstract}




\end{opening}
\section{Introduction}
Consider $GI/M/1/n$ queueing system denoting by $A(x)$ the probability
distribution function of interarrival time and by $\lambda$ the reciprocal
of the expected interarrival time,
$\alpha(s)=\int_0^\infty
\mbox{e}^{-sx}\mbox{d}A(x)$. The parameter of the service time distribution
will be denoted $\mu$, and load of the system is $\rho=\lambda/\mu$.
The size of buffer $n$ includes the position for server.
Denote also $\rho_m=\mu^m\int_0^\infty x^m\mbox{d}A(x)$, $m=1,2,...$,~
$(\rho_1=\rho^{-1})$.

The explicit representation for the loss probability
in terms of generating function was obtained by Miyazawa [12].
Namely, he showed that whenever the value of load $\rho$, the loss
probability $p_n$ always exists and has the representation
$$
p_n=\frac{1}{\sum_{j=0}^n\pi_j},
\eqno (1.1)
$$
where the generating function $\Pi(z)$ of $\pi_j$, $j=0,1,...$,
is the following
$$
\Pi(z)=\sum_{j=0}^\infty\pi_jz^j
=\frac{(1-z)\alpha(\mu-\mu z)}{\alpha(\mu-\mu z)-z},~~|z|<\sigma,
\eqno (1.2)
$$
$\sigma$ is the minimum nonnegative solution of the functional
equation $z=\alpha(\mu-\mu z)$.
This solution is the following. It belongs to the open interval
(0,1) if $\lambda<\mu$, and it is equal to 1 otherwise.

In the recent papers Choi and Kim [8] and Choi et al [9] study the questions
related to the asymptotic behavior of the sequence $\{\pi_j\}$ as
$j\to\infty$. Namely, they study asymptotic behavior of the loss probability
$p_n$, $n\to\infty$, as well as obtain the convergence rate of the stationary
distributions of the $GI/M/1/n$ queueing system to those of the
$GI/M/1$ queueing system as $n\to\infty$. The analysis of
[8] and [9] is based on the theory of analytic functions.

The approach of this paper is based on Tauberian theorems with remainder
permitting us to simplify the proof of the results of the mentioned
paper of Choi et al [9] as well as to obtain some new results
on asymptotic behavior of the loss probability.

For the asymptotic behavior of the loss probability in $M/GI/1/n$ queue
see Abramov [1], [2], Asmussen [6], Takagi [16],
Tomko [17], Willmot [18] etc.
For the asymptotic analysis of
more general than $M/GI/1/n$ queueing systems see Abramov [3],
Baiocchi [7] etc.

Study of the loss probability and its asymptotic analysis is motivated
by growing development of communication systems. The results of our
study can be applied to the problems of flow control, performance evaluation,
redundancy. For application of the loss probability to such kind of problems
see Ait-Hellal et al [4], Altman and Jean-Marie [5], Cidon et al
[10], Gurewitz et al [11].

\section{Auxiliary results. Tauberian theorems}
In this section we represent the asymptotic results of Takacs [15, p.22-23]
(see Lemma 2.1 below), and Tauberian theorems
of Postnikov [13, Section 25],
(see Lemmas 2.2 and 2.3 below).\\

Let $Q_j$, $j=0,1,...$, be a sequence of real numbers satisfying the
recurrent relation
$$
Q_n=\sum_{j=0}^n r_jQ_{n-j+1},
\eqno (2.1)
$$
where $r_j$, $j=0,1,...$, are nonnegative numbers, $r_0>0$, $r_0+r_1+...=1$, and
$Q_0>0$ is an arbitrary real number.

Denote $r(z)=\sum_{j=0}^\infty r_jz^j$, $|z|\le 1$, $\gamma_m=r^{(m)}(1-0)=\lim_{z\uparrow 1}
r^{(m)}(z)$,
where $r^{(m)}(z)$ is the $m$th derivative of $r(z)$. Then for $Q(z)=
\sum_{j=0}^\infty Q_jz^j$, the generating function of $Q_j$, $j=0,1,...$,
we have the following representation
$$
Q(z)=\frac{Q_0r(z)}{r(z)-z}.
\eqno (2.2)
$$
The statements below are known theorems on asymptotic behavior of the
sequence $\{Q_j\}$ as $j\to\infty$. Lemma 2.1 below joins two results by Takacs [15]: Theorem 5
on p. 22 and relation (35) on p. 23.\\

{\bf Lemma 2.1} (Takacs [15]). {\it If $\gamma_1<1$ then
$$
\lim_{n\to\infty}Q_n=\frac{Q_0}{1-\gamma_1}.
$$
If $\gamma_1=1$ and $\gamma_2<\infty$ then
$$
\lim_{n\to\infty}\frac{Q_n}{n}=\frac{2Q_0}{\gamma_2}.
$$
If $\gamma_1>1$ then
$$
\lim_{n\to\infty}\Big(Q_n-\frac{Q_0}{\delta^n[1-r'(\delta)]}\Big)=
\frac{Q_0}{1-\gamma_1},
$$
where $\delta$ is the least (absolute) root of equation $z=r(z)$}.\\

{\bf Lemma 2.2} (Postnikov [13]). {\it Let $\gamma_1=1$ and $\gamma_3<\infty$.
Then as $n\to\infty$
$$
Q_n=\frac{2Q_0}{\gamma_2}n+O(\log n).
$$}

{\bf Lemma 2.3} (Postnikov [13]). {\it Let $\gamma_1=1$, $\gamma_2<\infty$ and
$r_0+r_1<1$. Then as $n\to\infty$
$$
Q_{n+1}-Q_n=\frac{2Q_0}{\gamma_2}+o(1).
$$}

\section{The main results on asymptotic behavior of the loss probability}
Let us study (1.1) and (1.2) more carefully. Represent (1.2) as the difference of
two terms
$$
\Pi(z)=
\frac{(1-z)\alpha(\mu-\mu z)}{\alpha(\mu-\mu z)-z}=
\frac{\alpha(\mu-\mu z)}{\alpha(\mu-\mu z)-z}-
z\frac{\alpha(\mu-\mu z)}{\alpha(\mu-\mu z)-z}
$$
$$
=\widetilde\Pi(z)-z\widetilde\Pi(z),
\eqno (3.1)
$$
where
$$
\widetilde\Pi(z)=\sum_{j=0}^\infty\widetilde\pi_jz^j
=\frac{\alpha(\mu-\mu z)}{\alpha(\mu-\mu z)-z}.
\eqno (3.2)
$$
Note also that
$$
\pi_0=\widetilde\pi_0=1,
$$
$$
\pi_{j+1}=\widetilde\pi_{j+1}-\widetilde\pi_j,~~~j\ge 0.
\eqno (3.3)
$$
Therefore,
$$
\sum_{j=0}^n\pi_j=\widetilde\pi_n,
$$
and
$$
p_n=\frac{1}{\widetilde\pi_n}.
\eqno(3.4)
$$
Now the application of Lemma 2.1 yields the following\\

{\bf Theorem 3.1.} {\it In the case where $\rho<1$ as $n\to\infty$ we have
$$
p_n=\frac{(1-\rho)[1+\mu\alpha'(\mu-\mu\sigma)]\sigma^n}
{1-\rho-\rho[1+\mu\alpha'(\mu-\mu\sigma)]\sigma^n}+o(\sigma^{2n}).
\eqno (3.5)
$$
In the case where $\rho_2<\infty$ and $\rho=1$ we have
$$
\lim_{n\to\infty} np_n=\frac{\rho_2}{2}.
\eqno (3.6)
$$
In the case where $\rho>1$ we have
$$
\lim_{n\to\infty} p_n=\frac{\rho-1}{\rho}.
\eqno (3.7)
$$}

{\bf Proof.} Indeed, it follows from (3.1), (3.2) and (3.3) that $\widetilde \pi_0=1$, and
$$
\widetilde\pi_k=\sum_{i=0}^k\frac{(-\mu)^i}{i!}\alpha^{(i)}(\mu)\widetilde\pi_{k-i+1},
\eqno (3.8)
$$
where $\alpha^{(i)}(\mu)$ denotes the $i$th derivative of $\alpha(\mu)$. Note also that $\alpha(\mu)>0$,
the terms $(-\mu)^i\alpha^{(i)}(\mu)/i!$ are nonnegative for all $i\ge 1$, and
$$
\sum_{i=0}^\infty\frac{(-\mu)^i}{i!}\alpha^{(i)}(\mu)=\sum_{i=0}^\infty\int_0^\infty\mbox{e}^{-\mu x}\frac{(\mu x)^i}{i!}
\mbox{d}A(x)
$$
$$
=\int_0^\infty\sum_{i=0}^\infty\mbox{e}^{-\mu x}\frac{(\mu x)^i}{i!}\mbox{d}A(x)=1.
\eqno (3.9)
$$
Therefore one can apply Lemma 2.1. Then in the case of $\rho<1$ one can write
$$
\lim_{n\to\infty}\Big(\widetilde\pi_n-\frac{1}{\sigma^n[1+\mu\alpha'(\mu-\mu\sigma)]}\Big)=\frac{\rho}{\rho-1},
\eqno (3.10)
$$
and for large $n$ relation (3.10) can be rewritten in the form of the estimation
$$
\widetilde\pi_n=\Big[\frac{1}{\sigma^n[1+\mu\alpha'(\mu-\mu\sigma)]}+\frac{\rho}{\rho-1}\Big][1+o(\sigma^n)].
\eqno (3.11)
$$
In turn, from (3.11) for large $n$ we obtain
$$
p_n=\frac{1}{\widetilde\pi_n}=\frac{(1-\rho)[1+\mu\alpha'(\mu-\mu\sigma)]\sigma^n}
{1-\rho-\rho[1+\mu\alpha'(\mu-\mu\sigma)]\sigma^n}[1+o(\sigma^n)]
$$
$$
=\frac{(1-\rho)[1+\mu\alpha'(\mu-\mu\sigma)]\sigma^n}
{1-\rho-\rho[1+\mu\alpha'(\mu-\mu\sigma)]\sigma^n}+o(\sigma^{2n}).
$$
Thus (3.5) is proved. The limiting relations (3.6) and (3.7) follow immediately by application of Lemma 2.1.
Theorem 3.1 is proved.\\

The following two theorems improve limiting relation (3.6). From
Lemma 2.2 we have the following\\

{\bf Theorem 3.2.} {\it Assume that $\rho=1$ and $\rho_3<\infty$. Then as $n\to\infty$
$$
p_n=\frac{\rho_2}{2n}+O\Big(\frac{\log n}{n^2}\Big).
\eqno (3.12)
$$}

{\bf Proof.} The result follows immediately by application of Lemma 2.2.\\

Subsequently, from Lemma 2.3 we have\\

{\bf Theorem 3.3.} {\it Assume that $\rho=1$ and $\rho_2<\infty$. Then as $n\to\infty$
$$
\frac{1}{p_{n+1}}-\frac{1}{p_n}=\frac{2}{\rho_2}+o(1).
\eqno (3.13)
$$}

{\bf Proof}. The theorem will be proved if we show that for all $\mu>0$
$$
\alpha(\mu)-\mu\alpha'(\mu)<1.
\eqno (3.14)
$$
Taking into account (3.9) and the fact that $(-\mu)^i\alpha^{(i)}(\mu)/i!\ge 0$ for all $i\ge 0$,
one can write
$$
\alpha(\mu)-\mu\alpha'(\mu)\le 1.
\eqno (3.15)
$$
Thus, we have to show that for some $\mu_0>0$ the equality
$$
\alpha(\mu_0)-\mu_0\alpha'(\mu_0)=1
\eqno (3.16)
$$
is not a case. Indeed, since $\alpha(\mu)-\mu\alpha'(\mu)$ is an analytic function then,
according to the theorem on maximum absolute value of analytic function, the equality
$\alpha(\mu)-\mu\alpha'(\mu)=1$ is valid for all $\mu>0$. This means that
(3.16) is valid if and only if $\alpha^{(i)}(\mu)=0$ for all $i\ge 2$ and for all $\mu>0$,
and therefore $\alpha(\mu)$
is a linear function,
i.e. $\alpha(\mu)=c_0+c_1\mu$, where $c_0$ and $c_1$ are some constants.
However, since $|\alpha(\mu)|\le 1$ we obtain $c_0=1$, $c_1=0$.
This is a trivial case where the probability distribution function
$A(x)$ is concentrated in point 0. Therefore (3.16) is not a case, and hence (3.14) holds.
Theorem 3.3 is proved.\\

We have also the following\\

{\bf Theorem 3.4}. {\it Let
$\rho=1-\epsilon$, where $\epsilon>0$,
and $\epsilon n\to C>0$ as $n\to\infty$ and $\epsilon\to 0$.
Assume
that $\rho_3=\rho_3(n)$ is a bounded function and there exists $\widetilde\rho_2=\lim_{n\to\infty}\rho_2(n)$.
Then,
$$
p_n=\frac{\epsilon\mbox{e}^{-2C/\widetilde\rho_2}}{1-\mbox{e}^{-2C/\widetilde\rho_2}}[1+o(1)].
\eqno (3.17)
$$}

{\bf Proof.} It was shown in Subhankulov [14, p. 326] that if
$\rho^{-1}=1+\epsilon$, $\epsilon>0$ and $\epsilon\to 0$, $\rho_3(n)$ is a bounded function, and there exists
$\widetilde\rho_2=\lim_{n\to\infty}\rho_2(n)$
then
$$
\sigma=1-\frac{2\epsilon}{\widetilde\rho_2}+O(\epsilon^2),
\eqno (3.18)
$$
where $\sigma=\sigma(n)$ is the minimum root of the functional
equation $z-\alpha(\mu-\mu z)=0$, $|z|\le 1$, and where the
parameter $\mu$ and the function $\alpha(z)$, both or one of them,
are assumed to depend on $n$. Therefore, (3.18) is also valid
under the assumptions of the theorem. Then after some algebra one
can obtain
$$
[1+\mu\alpha'(\mu-\mu\sigma)]\sigma^n=
\epsilon\mbox{e}^{-2C/\widetilde\rho_2}[1+o(1)],
$$
and the result easily follows from estimation (3.11).\\

{\bf Theorem 3.5.} {\it Let
$\rho=1-\epsilon$, where $\epsilon>0$,
and $\epsilon n\to 0$ as $n\to\infty$ and $\epsilon\to 0$.
Assume
that $\rho_3=\rho_3(n)$ is a bounded function and there exists $\widetilde\rho_2=\lim_{n\to\infty}\rho_2(n)$.
Then
$$
p_n=\frac{\widetilde\rho_2}{2n}+o\Big(\frac{1}{n}\Big).
\eqno (3.19)
$$}

{\bf Proof.} The proof follows by expanding of the main term of asymptotic relation (3.17)
for small $C$.

\section{Discussion}
We obtained a number of asymptotic results related to the loss probability
for the $GI/M/1/n$ queueing system by using Tauberian theorems with remainder. Asymptotic relations (3.6)
and (3.7) of Theorem 3.1 are the same as correspondent asymptotic relations of Theorem 3 of [9].
Asymptotic relation (3.5) of Theorem 3.1 improves correspondent asymptotic relation of Theorem 3 of [9],
however it
can be deduced from Theorem 3.1 of [8] and the second equation on p. 1016 of [8].
Under additional condition
$\rho_3<\infty$ the statement (3.12) of Theorem 3.2 is new. It improves the result of [9] under $\rho=1$:
the remainder obtained in Theorem 3.2 is $O(\log n/n^2)$ whereas under condition $\rho_2<\infty$ the remainder
obtained in Theorem 3 of [9] is $o(n^{-1})$.
Asymptotic relation (3.13) of Theorem 3.3
coincides with intermediate asymptotic relation on p. 441 of [9]. Theorems 3.4 and 3.5 are new.
They provide
asymptotic results where the load $\rho$ is close to 1.

\section*{Acknowledgement}
The author thanks the anonymous referees for a number of valuable comments.
\section*{References}
\indent

[1] V.M.Abramov, {\em Investigation of a Queueing System with Service Depending on Queue-Length},
(Donish, Dushanbe, 1991) (in Russian).

\smallskip

[2] V.M.Abramov, On a property of a refusals stream, J. Appl. Probab. 34 (1997) 800-805.

\smallskip

[3] V.M.Abramov, Asymptotic behavior of the number of lost packets, submitted for publication.

\smallskip

[4] O.Ait-Hellal, E.Altman, A.Jean-Marie and I.A.Kurkova, On loss probabilities in presence of redundant
packets and several traffic sources, Perform. Eval. 36-37 (1999) 485-518.

\smallskip

[5] E.Altman and A.Jean-Marie, Loss probabilities for messages with redundant packets feeding a finite
buffer, IEEE J. Select. Areas Commun. 16 (1998) 778-787.

\smallskip

[6] S.Asmussen, Equilibrium properties of the $M/G/1$ queue, Z. Wahrscheinlichkeitstheorie 58 (1981) 267-281.

\smallskip

[7] A.Baiocchi, Analysis of the loss probability of the $MAP/G/1/K$ queue, part I: Asymptotic theory,
Stochastic Models 10 (1994) 867-893.

\smallskip

[8] B.D.Choi and B.Kim, Sharp results on convergence rates for the distribution of $GI/M/1/K$ queues as
$K$ tends to infinity, J. Appl. Probab. 37 (2000) 1010-1019.

\smallskip

[9] B.D.Choi, B.Kim and I.-S.Wee, Asymptotic behavior of loss probability in $GI/M/1/K$ queue as $K$ tends
to infinity, Queueing Systems 36 (2000) 437-442.

\smallskip

[10] I.Cidon, A.Khamisy and M.Sidi, Analysis of packet loss processes in high-speed networks, IEEE Trans.
Inform. Theory 39 (1993) 98-108.

\smallskip

[11] O.Gurewitz, M.Sidi and I.Cidon, The ballot theorem strikes again: Packet loss process distribution,
IEEE Trans. Inform. Theory 46 (2000) 2588-2595.

\smallskip

[12] M.Miyazawa, ~~Complementary~ generating~ functions for the $M^X/GI/1/k$ and $GI/M^Y/1/k$ queues and
their application to the comparison for loss probabilities, J. Appl. Probab. 27 (1990) 684-692.

\smallskip

[13] A.G.Postnikov,~~ Tauberian~ theory and its application. Proc. Steklov Math. Inst. 144 (1979) 1-148.

\smallskip

[14] M.A.Subhankulov, ~~{\em Tauberian~ Theorems with Remainder}. (Nauka, Moscow, 1976) (in Russian).

\smallskip

[15] L.Takacs, {\em Combinatorial Methods in the Theory of Stochastic Processes}. (Wiley, New York, 1967).

\smallskip

[16] H.Takagi, {\em Queueing Analysis}, Vol. 2. (Elsevier Science, Amsterdam, 1993).

\smallskip

[17] J.Tomko, One limit theorem in queueing problem as input rate increases infinitely. Studia Sci. Math.
Hungarica 2 (1967) 447-454 (in Russian).

\smallskip

[18] G.E.Willmot, A note on the equilibrium $M/G/1$ queue-length, J. Appl. Probab. 25 (1988) 228-231.

\end{article}
\end{document}